\documentclass[12pt,a4paper,leqno]{amsart}

\usepackage{amssymb, amscd}

\def\Hom{\operatorname{Hom}}

\font\cyr=wncyr10 scaled \magstep1%
\def\Sh{\text{\cyr Sh}}

\newtheorem{thm}{Theorem}[section]
\newtheorem{conj}[thm]{Conjecture}
\newtheorem*{thm*}{Theorem}

\newtheorem*{remark*}{Remarks}

\newtheorem*{defn*}{Definition}
\newtheorem*{claim*}{Claim}

\theoremstyle{definition}
\newtheorem{remark}[thm]{Remark}

\newtheorem{quest}[thm]{Question}

 \renewcommand{\sectionmark}[1]{}

\newcommand{\ve}{\varepsilon}

\newcommand{\iy}{\infty}

\newcommand{\Dl}{\Delta}

\newcommand{\ov}{\overline}

\begin{document}

\title[Brauer--Siegel theorem for elliptic surfaces]
{Brauer--Siegel theorem for elliptic surfaces}
\author[B. \`E. Kunyavski\u i]{B. \`E. Kunyavski\u i}
\address{Kunyavski\u\i : Department of Mathematics,
Bar-Ilan University, Ramat Gan, Israel}
\email{kunyav@macs.biu.ac.il}
\author[M. A. Tsfasman]{M. A. Tsfasman}
\address{Tsfasman: French-Russian Poncelet Laboratory;
Institut de Math\'e\-ma\-tiques de Luminy; Independent University
of Moscow; and Institute for Information Transmission Problems}
\email{tsfasman@iml.univ-mrs.fr}

%\email{}

\thanks{This research was supported in part by the
French-Israeli grant 3-1354 and the Russian-Israeli grant RFBR
06-01-72004-MSTIa. The first named author was also supported in
part by the Minerva Foundation through the Emmy Noether Institute
of Mathematics. The second named author was also supported in part
by RFBR 02-01-22005, 02-01-01041, 06-01-72550-CNRSa, 07-01-00051,
and INTAS 05-96-4634.}

\begin{abstract}
We consider higher-dimensional analogues of the classical
Brauer-Siegel theorem focusing on the case of abelian varieties
over global function fields. We prove such an analogue in the case
of constant families of elliptic curves and abelian varieties.
\end{abstract}
%\subjclass{30C80, 30D40}

\maketitle

\begin{center}
{\em {To our teachers V.E.~Voskresenski\u\i \ and Yu.I.~Manin
\\[2pt]

to their 80th and 70th birthdays, respectively}}

\end{center}

\baselineskip 20pt

\section{Introduction} \label{sec:intro}

The classical Brauer--Siegel theorem, which is one of the
milestones of the number theory of the past century, reflects deep
connections between algebraic, arithmetical, analytic, and (in the
function field case) geometric properties of global fields. Not
only is the theorem a working tool in a variety of problems
concerning number and function fields, but the underlying ideas
have been recently put into much broader context expanding far
beyond number theory (see, for example, \cite{ST}).

Recall that the theorem describes the asymptotic behaviour of the
product of two important arithmetic invariants of a number field
$K$, the class number $h(K)$ and the regulator $R(K)$, as the
discriminant $d(K)$ tends to infinity. More precisely, it says
that the ratio $r=\log(hR)/\log(\sqrt{|d|})$ tends to 1 provided
at least one of the following conditions is satisfied: 1) the
degree $n=[K:\mathbb Q]$ remains the same for all $K$'s in the
sequence of fields under consideration; 2) $n/\log(|d|)$ tends to
0 and all $K$'s are normal. Even in this not-so-effective form
there are many useful applications. Some effective versions of the
theorem are known in several particular cases (see \cite{St} and
references therein).

A natural question whether the statement of the theorem still
holds when none of conditions 1) and 2) is satisfied, or under
some weaker assumptions, remained widely open until recently. In
the paper \cite{TV2} there were obtained some asymptotic bounds on
$r$ generalizing the statement of the Brauer--Siegel theorem.
These techniques, together with those of an earlier paper
\cite{TV1} led to a new concept of infinite global field which is
an important object for further investigation. Combined with
Weil's ``explicit formulae'' (see \cite{LT}), they yielded quite a
few concrete arithmetic applications, like new estimates for
regulators. Note that even more general approach was used in a
recent paper \cite{Zy} where the normality assumption on $K$ was
weakened.

The above mentioned results present the state of the art in the
research area concentrated around the classical Brauer--Siegel
theorem. In the present paper we make an attempt to treat some new
problems arising from these achievements. Namely, one can think
about higher dimensional analogues of the Brauer--Siegel theorem.
In particular, if $E$ is a commutative algebraic group defined
over a global field $K$, one can define an analogue of the class
number $h(E)$ and the regulator $R(E)$. Moreover, the classical
analytical class number formula of Dirichlet admits higher
dimensional analogues both for algebraic tori \cite{Shyr} and,
conjecturally, for abelian varieties (Birch and Swinnerton-Dyer).
This motivates the study of asymptotic behaviour of $h(E)R(E)$ in
appropriately chosen families of groups $E$ when the
``discriminant'' $d(E)$ tends to infinity. In the case where $E$
is an abelian variety, recent work of Hindry and Pacheco contains
quite a new approach to this kind of asymptotic problems, both in
the number field case \cite{Hi} and in the function field case
\cite{HP}. This work was an additional motivation for publishing
our results because the approach of Hindry and Pacheco is, in a
sense, ``orthogonal'' to ours: loosely speaking, they consider
``vertical'' families of abelian varieties (say, in the function
field case the genus of the underlying curve $X$ is fixed and the
conductor of the abelian variety grows) while we consider
``horizontal'' families where the genus of $X$ tends to infinity.

\section{Main theorem} \label{sec:main}

We fix the ground field $k={\mathbb F}_q$ and consider a (smooth,
projective, geometrically irreducible) curve $X/\mathbb F_q$ of
genus  $g$.
%($g\to\iy$).
Let $K=\mathbb F_q(X)$, and let $E/K$ be a (smooth, connected)
commutative algebraic $K$-group.  Our goal is to study asymptotic
behaviour of the ``class number" $h(E)$ as $g\to\iy.$ In the
present paper we focus on the particular case where $E=A$ is an
abelian variety (see, however, Section \ref{sec:gen} for
the case where $E$ is an algebraic $K$-torus). Let
$\Sh:=\vert\Sh (A)\vert$ be the order of the Shafarevich--Tate
group of $A$, and $\Dl$ the determinant of the Mordell--Weil
lattice of $A$ (cf. \cite{Mi}, \cite{Hi}). In this section we
consider the most trivial ``constant" case, i.e. $E\cong
E_0\times_{\mathbb F_q}K$  where $E_0$ is an $\mathbb F_q$-group;
see Section \ref{sec:gen} for a more general setting.

To state our main result, we recall some notation from \cite{TV1}.
If the ground curve $X=X_0$ varies in a family $\{X_i\}$, we
denote by $g_i$ the genus of $X_i$ ($g_i\to\iy$), by $N_m(X_i)$
the number of ${\mathbb F}_{q^m}$-points of $X_i$, and we always
assume that for every $m\ge 1$ there exists a limit
$\beta_m:=\lim_{i\to\iy}\frac{N_m(X_i)}{g_i}$. Such families are
called {\it asymptotically exact}; any family contains an
asympotically exact subfamily; any tower (i.e. a family such that
$k(X_i)\subset k(X_{i+1})$ for every $i$) is asymptotically exact;
see \cite{Ts}, \cite{TV1} for more details. We shall often drop
the index $i$ if this does not lead to confusion.

\begin{thm}\label{th:main}
Let $E = E_0\times_{{\mathbb F}_q}K$ where $E_0$ a fixed elliptic
${\mathbb F}_q$-curve.
%%%%%%%%%%%%%%%%%%%%%%%%%
Let $K$ vary in an asymptotically exact family, and let $\beta_m$
be the corresponding limits.
%%%%%%%%%%%%%%%%%%%%%%%%%%%%%%%%%%%%
Then
$$\lim_{i\to\iy}\frac{1}{g_i}\log_q(\Sh\cdot\Dl) = 1-\sum_{m=1}^\iy\beta_m\log_q\frac{N_m(E_0)}{q^m},$$
where $N_m(E_0)=|E_0(\mathbb F_{q^m})|$.
\end{thm}

\begin{proof}
%First consider the case where $A_0=E_0$ is an elliptic curve over
%$k,$\ $A=E=E_0\times_k K$ is a constant elliptic curve over
%$K=k(X).$

Denote by $\omega_j$\ $(j=1,\dots,2g)$ the eigenvalues of
Frobenius acting on $H^1(X),$ and by $\psi_1,\psi_2$ the
eigenvalues of Frobenius acting on $H^1(E_0).$ We have
$\omega_j\ov \omega_j=\psi_1\psi_2=q.$

Put $t=q^{-s}$ and consider the Hasse--Weil  $L$-function of
$E/K.$ According to the Birch and Swinnerton-Dyer conjecture
(which, under our hypotheses, is a theorem \cite{Mi}, \cite{Oe}),
the value of $L_E(t)/(1-qt)^r$  at $t=q^{-1}$ equals
$q^{1-g}\cdot\Sh\cdot\Dl/[\#E_0(k)]^2.$  Here $r$ is the rank of
$E(K)/E(K)_{\text{tors}}$; this number is equal to the number of
pairs $(i,j)$ such that $\psi_i=\omega_j$ (loc. cit.). This gives
us Milne's formula
$$\Sh\cdot \Dl=q^g\prod_{\omega_j\ne\psi_i}\left(1-
\frac{\psi_i}{\omega_j}\right).$$ It is convenient to put
$\psi_i=\alpha_i\sqrt q,$\ $\omega_j=\gamma_j\sqrt q,$  and,
taking into account that the Frobenius roots can be written as
conjugate pairs, to write the above formula as
\begin{equation}\label{1}
\Sh\cdot \Dl=q^g\prod_{\alpha_i\ne
1/\gamma_j}(1-\alpha_i\gamma_j).
\end{equation}

Set $\alpha_1=\alpha,$\ $\alpha_2=\ov\alpha.$ First consider the
case where $r=0.$ Then the right-hand side of \eqref{1} can be
written as $q^gP_X(\alpha/\sqrt q)P_X(\ov\alpha/\sqrt q),$ where
$P_X(t)$ is the numerator of the zeta-function of $X$ :
$$Z_X (t)=\frac{P_X(t)}{(1-t)(1-qt)}.$$
Hence the right-hand side of \eqref{1} equals
$$q^g\left[\left(1-\frac{\alpha}{\sqrt q}\right)(1-\alpha\sqrt q)Z_X\left(\frac{\alpha}{\sqrt
q}\right)\left(1-\frac{\ov\alpha}{\sqrt q}\right)(1-\ov\alpha\sqrt
q)Z_X\left(\frac{\ov\alpha}{\sqrt q}\right)\right].$$ We now write
$Z_X(t)=\prod\limits_{m=1}^\iy(1-t^m)^{-B_m},$ then we have
$\beta_m=\lim\limits_{g\to\iy}\frac{B_m}{g}$ (by our assumption,
the limit exists), and we get
\begin{align*}
\lim_{g\to\iy}\frac{1}{g}\log_q(\Sh\cdot\Dl)&=1+\log_q\left(\prod_{m=1}^\iy\left(1-\frac{\alpha^m}{
q^{\frac{m}{2}}}\right)^{-\beta_m}\left(1-\frac{\ov\alpha^m}{q^{\frac{m}{2}}}\right)^{-\beta_m}\right)\\
&=1-\sum_{m=1}^\iy\beta_m\log_q\left(1+\frac{1}{q^m}-\frac{\alpha^m+\ov\alpha^m}{q^{\frac{m}{2}}}
\right)\\ &=1-\sum_{m=1}^\iy\beta_m\log_q\frac{N_m}{q^m}
\end{align*}
(here $N_m=|E_0(\mathbb F_{q^m})|$, and the last equality follows from the Weil formula). Note that the
series on the right-hand side converges according to \cite{Ts}.
Indeed, we know that the series
$\sum\limits_{m=1}^\iy\frac{m\beta_m}{q^{\frac{m}{2}}-1}$
converges \cite[Cor.1]{Ts}. We have
$\frac{N_m}{q^m}=1+q^{-m}-\frac{\alpha^m+\ov\alpha^m}{q^{\frac{m}{2}}}.$
Fix $m>m_0$ big enough. Put
$x=\frac{\alpha^m+\ov\alpha^m}{q^{\frac{m}{2}}}-q^{-m}.$ Since
$|\alpha^m+\ov\alpha^m|\le  2,$ we have
$$\left|\log_q\frac{N_m}{q^m}\right|=|\log_q(1-x)|\le
c\sum_{n=1}^\iy\left(q^{-\frac{m}{2}}\right)^n\le c'q^{-\frac{m}{2}}\le
c'\frac{m}{q^{\frac{m}{2}}-1}.$$
Hence the series $\sum\beta_m\log_q\frac{N_m}{q^m}$ converges.

Let us now consider the case where $r>0.$ Our key observation is
that as $g\to\iy,$ the rank cannot grow as fast as $g,$ i.e., we
always have $\lim\limits_{g\to\iy}\frac{r}{g}=0.$

Indeed, if $\lim\limits_{g\to\iy}\frac{r}{g}=c>0, $ then there is
at least one multiple Frobenius root $\omega_j=\psi_1$ or $\psi_2$
with multiplicity $\ge cg.$ Hence the Weil measure (cf.
\cite{TV1})
$$
\mu_\Omega=\frac{1}{g}\sum_{j=1}^{2g}\delta_{\gamma_j}\quad
(\text{where } \delta_{\gamma_j} \text{is the Dirac measure)}
$$
tends (as $g\to\iy)$ to a measure that is greater than or equal to
$c\delta_{\gamma_j}.$ But according to \cite[Th.2.1]{TV1}, the
limit measure $\mu=\lim\limits_{g\to\iy}\mu_\Omega$ must have a
continuous density, contradiction.

%%%%%%%%%%%%%%%%%%%%%%%%%%%%%%%%%%%%%%%
(As pointed out by the referee, this observation might happen to
be deducible from the ``explicit formulae'' for elliptic curves, see,
e.g., \cite{Br}.)
%%%%%%%%%%%%%%%%%%%%%%%%%%%%%%%%%%%%%%

Thus, in the general case where $r>0, $ we get the required result
as follows.

Let us introduce an auxiliary function $\delta (g)=1+\ve(g)$ such
that $\lim\limits_{g\to\iy}\ve(g)=0$ and
$\lim\limits_{g\to\iy}\left(\frac{r\log\ve(g)}{g}\right)=0.$
%ETO NEVERNO! (Since $\frac{r}{g}\to0,$ we can take,
%say, $\ve(g)=g^{-1}.)$
Let
$$F(g)=q^gP_X(\delta(g)\alpha/\sqrt
q)P_X(\delta(g)\ov\alpha/\sqrt q).$$ We have, on the one hand,
\begin{equation}\label{2}
\lim\limits_{g\to\iy}\frac{\log
_qF(g)}{g}=1-\sum_{m=1}^\iy\beta_m\log_q\left(\frac{N_m}{q^m}\right),
\end{equation}
and, on the other hand,
$$
\lim_{g\to\iy}\frac{1}{g}\log_qF(g)=\lim_{g\to\iy}\left(\frac{1}{g}\log_q(\Sh\cdot\Delta)\right).$$
To prove the last equality, we write {\scriptsize{
\begin{align*}
&F(g)=q^g\prod_{j=1}^{2g}(1-\alpha\gamma_j\delta(g))(1-\ov\alpha\gamma_j\delta(g))=\delta(g)^{4g}q^g\prod_{
j=1}^{2g}\left(\frac{1}{\delta(g)}-\alpha\gamma_j\right)\left(\frac{1}{\delta(g)}-\ov\alpha\gamma_j\right)\\
&=\delta(g)^{4g}q^g\prod_{\gamma_j=1/\alpha}\left(\frac{1}{\delta(g)}-\alpha\gamma_j\right)
\left(\frac{1}{\delta(g)}-\ov\alpha\gamma_j\right)\cdot\prod_{\gamma_j\ne
1/\alpha}\left(\frac{1}{\delta(g)}-\alpha\gamma_j\right)
\left(\frac{1}{\delta(g)}-\overline\alpha\gamma_j\right)\\
&=\delta(g)^{4g}\left(\frac{1}{\delta(g)}-1\right)^r\cdot
q^g\cdot\prod_{\gamma_j\ne 1/\alpha}(1-\alpha
\gamma_j\delta(g))(1-\ov\alpha\gamma_j\delta(g))\cdot
\frac{1}{\delta(g)^{4g-r}}\\
&=(1-\delta(g))^r\cdot q^g\cdot\prod_{\gamma_j\ne
1/\alpha}(1-\alpha\gamma_j\delta(g))(1-\ov\alpha\gamma_j\delta(g)).
\end{align*}
}} Hence
\begin{align*}
\lim_{g\to\iy}\frac{1}{g}\log_qF(g)&=\lim_{g\to\iy}\left(\frac{1}{g}\log_q(1-\delta(g))^r\right)+\lim_{
g\to\iy}\left(\frac{1}{g}\log_q(\Sh\cdot\Delta)\right)\\
&=\lim_{g\to\iy}\frac{r\log_q\ve(g)}{g}+\lim_{g\to\iy}\left(\frac{1}{g}\log_q(\Sh\cdot\Delta)\right)\\&=\lim
_{g\to\iy}\left(\frac{1}{g}\log_q(\Sh\cdot\Delta)\right).
\end{align*}

Note that the series
$$
\sum_{m=1}^\iy\beta_m\log_q\left(1+\frac{1}{q^m}+\frac{(\delta(g)\alpha)^m+(\delta(g)\ov\alpha)^m}
{q^{\frac{m}{2}}}\right)
$$
converges for every fixed $\delta(g)$ sufficiently close to 1. Hence the passage
to the limit in \eqref{2} is legitimate.
\end{proof}

A direct analogue of Theorem \ref{th:main} is true for constant
abelian varieties of arbitrary dimension.

\begin{thm}\label{th:ab}
Let $A = A_0\times_{{\mathbb F}_q}K$ where $A_0$ a fixed abelian
${\mathbb F}_q$-variety of dimension $d$.
%%%%%%%%%%%%%%%%%%%%%%%%%
Let $K$ vary in an asymptotically exact family, and let $\beta_m$
be the corresponding limits.
%%%%%%%%%%%%%%%%%%%%%%%%%%%%%%%%%%%%
Then
$$\lim_{i\to\iy}\frac{1}{dg_i}\log_q(\Sh\cdot\Dl) = 1-\sum_{m=1}^\iy\beta_m\log_q\frac{{N_m(A_0)}^{1/d}}{q^m},$$
where $N_m(A_0)=|A_0(\mathbb F_{q^m})|$.
\end{thm}

%Finally, let us consider the case where $d=\dim A>1.$
%%%%%%%%%%%%%%%%%%%%%%%%%%%%%%%%%%%%%%%%%%%%%%%%%%%%%%%%%%%%%%%%%%%%%
\begin{proof}

The proof goes as for elliptic curves,
{\it mutatis mutandis}. The value of $L_A(t)/(1-qt)^r$ at $t=q^{-1}$
equals $q^{1-dg}\cdot\Sh \cdot \Delta /(\#A_0(k)\cdot\#A_0^{\vee}(k))$, where
$A_0^{\vee}$ stands for the dual abelian variety.
According to \cite[Th.~3]{Mi}, this leads to
%%%%%%%%%%%%%%%%%%%%%%%%%%%%%%%%%%%%%%%%%%%%%%%%%%%%%%%%%%%%%%%%%%%
a formula similar to~\eqref{1}
$$\Sh\cdot\Delta=q^{dg}\prod_{\alpha_i\ne 1/\gamma_j}(1-\alpha_i\gamma_j),$$
where $\alpha_i$\ $(i=1,\dots,2d)$ are the (normalized) Frobenius
roots of $A_0.$ Therefore the case $r_A=0$ is treated, word for word, as in
the case $d=1.$ If $r_A>0,$ we have to prove that
$\frac{r_A}{g_X}\to 0$ as $g_X\to\iy$, and then apply the same
argument as for elliptic curves. Assume the contrary, i.e.,
$\lim\limits_{g_X\to\iy}\frac{r_A}{g_X}=c>0.$ Note that the
Mordell--Weil group $A(K)/A(K)_{\text{tors}}$ is isomorphic to
$\Hom_k(J_X,A_0).$ This implies that at least one Frobenius root
of $J_X$ (or of $X$, which is the same) appears with the
multiplicity proportional to $g.$ As in the one-dimensional case,
we then consider the Weil measure $\mu_\Omega$ and see that its
limit as $g\to\iy$ has discontinuous density which contradicts
\cite{TV1}.

The theorem is proved.
\end{proof}

\section{Generalizations} \label{sec:gen}

In this section we shall describe some possible generalizations of
Theorem \ref{th:main}. To make our approach more clear, we shall
first restrict ourselves to considering the case where $E$ is an
elliptic $K$-curve. Denote by $\mathcal E$ the corresponding
elliptic surface (this means that there is a proper connected
smooth morphism $f\colon \mathcal E \to X$ with the generic fibre
$E$). Assume that $f$ fits into an infinite Galois tower, i.e.
into a commutative diagram of the following form:
\begin{equation}
\CD \mathcal E=\mathcal E_0   @<<<  \mathcal
E_1  @<<< \dots @<<< \mathcal E_j @<<< \dots \\
@VVfV  @VVV & &   @VVV \\
X=X_0   @<<<  X_1  @<<<  \dots  @<<<  X_j @<<< \dots ,
\endCD
\label{eq:CD}
\end{equation}
where each lower horizontal arrow is a Galois covering. Let us
introduce some notation. For every $v\in X$, let $E_v=f^{-1}(v)$,
let $r_{v,i}$ denote the number of points of $X_i$ lying above
$v$, $\beta_v=\lim_{i\to\iy}r_{v,i}/g_i$ (we suppose the limits
exist). Furthermore, denote by $f_{v,i}$ the residue degree of a
point of $X_i$ lying above $v$ (the tower being Galois, this does
not depend on the point), and let $f_v=\lim_{i\to\iy}f_{v,i}$. If
$f_v=\iy$, we have $\beta_v=0$. If $f_v$ is finite, denote by
$N(E_v,f_v)$ the number of ${\mathbb F}_{q^{f_v}}$-points of
$E_v$. Finally, let $\tau$ denote the ``fudge'' factor in the
Birch and Swinnerton-Dyer conjecture (see \cite{Ta} for its
precise definition). Under this setting, we dare formulate the
following

\begin{conj} \label{conj:Galois}
Assuming the Birch and Swinnerton-Dyer conjecture for elliptic
curves over function fields, we have
$$\lim_{g\to\iy}\frac{1}{g}\log_q(\Sh\cdot\Dl\cdot\tau) = 1-\sum_{v\in X}\beta_v\log_q\frac{N(E_v,f_v)}{q^{f_v}}.$$
\end{conj}

\begin{remark} One can check that in the constant case
Conjecture \ref{conj:Galois} is consistent with Theorem
\ref{th:main}. The first nontrivial case to be considered is that
of an isotrivial elliptic surface.
\end{remark}

%\section{Further generalizations} \label{sec:gen}

Here are some questions for further investigation.

\begin{quest}
How can one formulate an analogue of Conjecture \ref{conj:Galois}
for more general towers when diagram (\ref{eq:CD}) does not
commute? for more general families when there are no upper
horizontal arrows in diagram (\ref{eq:CD})?
\end{quest}

With an eye towards even further generalizations of the
Brauer--Siegel theorem to arbitrary commutative algebraic groups,
the next extreme case to be considered is that of algebraic tori.
In that case the analogues of the class number and the regulator
are known \cite{Ono}, \cite{Vo}. Moreover, there is an analogue of
the analytic class number formula of Dirichlet established in
\cite{Shyr} for tori over number fields. Together with Theorem
\ref{th:ab}, this motivates the following

\begin{conj}\label{conj:tor} Let $T = T_0\times_{{\mathbb
F}_q}K$, where $T_0$ is a fixed ${\mathbb F}_q$-torus. Then
$$\lim_{g\to\iy}\frac{1}{g}\log h(T) = \lim_{g\to\iy}\frac{1}{g}\log\sqrt{\mathcal
D_T}-\sum_{m=1}^\iy\beta_m\log_q\frac{N_m(T_0)}{q^{md}},$$ where
$d=\dim T,$\ $N_m(T_0)=|T_0(\mathbb F_{q^m})|,$\ $\mathcal D_T$ is
the ``quasi-discriminant'' of $T$ $($cf. \cite{Shyr}$)$, and all
other notation is as in the previous sections.
\end{conj}

\bigskip

%NB: The setting of \cite{Shyr} is to be modified in the function
%field case (see \cite{Oe1}). The definition of the
%quasi-discriminant is to be rethought.

%Idea: to relate it to the Swan conductor of $T.$ One can expect
%that $\mathcal D(T)$ can be expressed in terms of $\mathcal D_k,$\
%$\mathcal D_L$ \ $(L=$ the minimal splitting field of $T$), and
%the Artin conductor of the character module $\hat T.$
%\end{comment}

%\medskip

\noindent {\it Acknowledgement.} A substantial part of this work
was done during the visits of the first named author to the
Mediterranean University and the Institute of Mathematics of
Luminy in 2003 and 2007 and the visit of the second named author
to Bar-Ilan University in 2005. Hospitality and support of these
institutions are gratefully appreciated. The authors thank A.~
Zykin for useful discussions and the referee for helpful remarks.

\end{document}